**PAPER • OPEN ACCESS**

# A strategic planning of a digital copy (an enterprise) as a task of control a dynamic system

View the article online for updates and enhancements.

# A strategic planning of a digital copy (an enterprise) as a task of control a dynamic system

S N Masaev[1,2,6], V V Vingert[1], A V Bogdanov[3] and Y K Salal[4,5]

[1] Siberian Federal University, pr. Svobodnyj, 79, Krasnoyarsk, 660041, Russia
[2] Control Systems LLC, 86 Pavlova Street, Krasnoyarsk, 660122, Russia
[3] The Main Directorate of EMERCOM of Russia for Krasnoyarsk Territory, Mira Ave., 68, Krasnoyarsk, 660122, Russia
[4] South Ural State University, pr. Lenina 76, Chelyabinsk, 454080, Russia
[5] University of Al-Qadisiyah, Diwaniyah, 58002, Iraq

[6] E-mail: faberi@list.ru

**Abstract**. The area of research includes control theory, dynamic systems, parameters of the external environment, mode, integral indicators, strategy. The general problem of assessing the state of large economic objects (enterprises) is revealed. There is no unified assessment a control of through strategic planning. The article proposes an integral indicator method for a unified assessment of the enterprise. The enterprise is formalized as a digital copy. The digital copy includes all business processes at any given time. The digital copy is presented as a multidimensional dynamic system. Dimension of multidimensional dynamic system is 1.2 parameters This allows you to estimate the modes operation of the enterprise in the normal mode and in the mode for control of strategic planning.

## 1. Introduction

Traditionally, strategic planning and goal planning are considered related processes, but this is not entirely true. Strategic planning and goal planning have a significant difference. Target control appeared in the 18th century in France, Germany and in the late 18th century in the United States. It appeared to control the activities of civil servants and the military. The planning of the strategy is commemorated in Ancient Greece, but it has been formalized for the control of an economic object since the middle of the 20th century. Since the middle of the 20th century, strategy planning and targets classifier have traditionally complemented each other for control an economic object (enterprise).

Strategic control originates from the pioneering work of A. Chandler in 1962. He indicated that the implementation of the economic strategy of an enterprise is possible only after the formation of goals and ways to achieve them [1]. Strategy development is inextricably linked to the following chain of control decisions: goal >> mission >> vision >> private >> strategic >> goals >> programs >> activities [1,2].

From the beginning of the 20th century to the beginning of the 21st century, the activities of economic objects (enterprises) have become much more complicated. For example, the well-known world leaders Coca-Cola, McDonald's, Boeing, BMW etc. have a huge system of business processes. Number of business processes is so large. It is more than all processes in a small state. Traditionally, the control of such large economic entities is carried out through various types of organizational structures with a constant increase in automation techniques. It contains a large data dimension and cannot be specified by traditional methods: intersectoral balances, vector, parametric and neural network modeling, agent-based approach, etc. [8-10] Simple control decisions become difficult to implement. For example, control a large enterprise through strategy planning. The problem is that the linking of the parameters of an economic object (enterprise) with existing strategies is not formalized.

There is no universal estimate of all dimensions of big data. Such an assessment can be an integral indicator. Therefore, a purpose of research is to estimate the state of an economic object (enterprise) as a multidimensional dynamic system in the basic mode of operation and its control mode, through the strategic planning with unknown parameters of the external environment.

## 2. Method an integral indicator

It is enough to imagine the economic activity of the enterprise (digital copy) as $S=\{T,X\}$, where $T=\{t:t=1,...,T_{max}\}$ - a lot of time points with a selected interval for analysis; $X$ - set of system parameters; $x(t)=[x^1(t),x^2(t),...,x^n(t)]^T \in X$ – $n$ – vector of indicators corresponding to the state of the system. Indicators of the vector $x^i(t)$ - the value of financial expenses and income of the enterprise. The dimension of the system $n$ is 1.2 million parameters. Based on the parameters $X$ and $T$, we consider our economic object (digital copy) a multidimensional dynamic system (hereinafter referred to as the system).

The analysis of the system at the moment $t$ is performed $x(t)$ on the basis $k$ of previous measures. The parameter $k$ is the length of the time series segment (accepted $k=6$ months in the work). Then we have a matrix

$$X_k(t) = \begin{bmatrix} x^T(t-1) \\ x^T(t-2) \\ \ldots \\ x^T(t-k) \end{bmatrix} = \begin{bmatrix} x^1(t-k) & x^2(t-k) & \cdots & x^n(t-k) \\ x^1(t-k) & x^2(t-k) & \cdots & x^n(t-k) \\ \ldots & \ldots & \ldots & \ldots \\ x^1(t-k) & x^2(t-k) & \cdots & x^n(t-k) \end{bmatrix} \quad (1)$$

$$R_k(t) = \frac{1}{k-1} \overset{o}{X_k^T}(t) \overset{o}{X_k}(t) = \left\| r_{ij}(t) \right\|, \quad (2)$$

$$r_{ij}(t) = \frac{1}{k-1} \sum_{l=1}^{k} \overset{o}{x^i}(t-l) \overset{o}{x^j}(t-l), \quad i,j=1,...,n, \quad (3)$$

where $t$ are the time instants, $r_{ij}(t)$ are the correlation coefficients of the variables $x^i(t)$ и $x^j(t)$ at the time instant $t$.

Next we form one of the four integral indicators – the sum of the absolute indicators of the correlation coefficients. It is indicator for express estimation of the correlation of system parameters $G_i(t)$:

$$R_i(t) = G_i(t) = \sum_{j=1}^{n} |r_{ij}(t)|. \quad (4)$$

The state of the entire system is calculated as:

$$G = \sum_{t=1}^{T=\max} \sum_{i=1}^{n} G_i(t). \quad (5)$$

The strategic planning $V$ is a set of system strategy of enterprise function (business-process), which can be represented $V_i^k$ as a set of strategic planning:

$$V = \sum_{t=1}^{T=\max} \sum_{i=1}^{n} V_i^k(t). \quad (6)$$

We will carry out the identification of the performed functions of the system with a strategy for the function (business-process). Each strategic planning $V_i^k$ is characterized by the business-process of the enterprise:

$$V_i^k = \sum_{i=1}^{n} v_i^j(x_j^i) \to \min, \qquad (7)$$

where $v_i^j$ is name of a strategy (compliance $x_j^i$ is $v_i^j$ set as 1-yes, 0-no); $x_j^i$ - the costs of the $i$ - the strategy for the function (business process) $j$. Control method is $V$. It is strategy $v_i^j$ so have a vector of strategic planning $v(t) = \left[v^1(t), v^2(t),...,v^n(t)\right]^T \in V$ - $n$ - dimension. Then $V = \sum_{t=1}^{T=\max} \sum_{i=1}^{n} V_i^k(t)$.

Payment of the functional duties of employees of the economic system is limited by resources C, then C(X)≤C. This restriction applies to all subsystems of the researched system.

The implementation of the method is performed in the author's complex of programs.

## 3. Characteristics of the research objects

Business processes of a wood processing enterprise in the Krasnoyarsk Territory are considered. The company operates in the Severo-Yenisei region and the city of Krasnoyarsk. The main business processes of the enterprise: 1) procurement of 800 thousand cubic meters of wood in the Severo-Yenisei region, 2) delivery along the Yenisei river to Krasnoyarsk and 3) production: floor boards, glued beams, calibrated boards, etc. After 1.5 years the company going to introduce a strategic control system and going to expand the volume of roundwood processing from 800 thousand to 1.5 million cubic meters.

At this enterprise, a control loop is being introduced through the strategic planning ($V_{strat\_plan}$): product strategy, corporate strategy, operational strategy, management (control) strategy, resource strategy, etc. [6]

## 4. Experiment result

Initial calculation data: $n$=1.2 million values, $X$=5,641,442 thousand rubles, control is set through the choice of strategies ($V_{strat\_plan}$). From the 7th period, incentive payments are added to specialists in all production areas to track the implementation of the strategy. Also, a business trip of the manager is paid for training the selected method of control at the enterprise (dynamic system). The calculation algorithm is 435 minutes.

A table 1 shows the experiment result of estimating the control mode $V_i(t)$ through the strategic planning.

**Table 1.** Regimes: $V_{(basic\ mode)}$ and $V_{(strat\_plan)}$.

| $t$ | $V_{(basic\ mode)}$ | $V_{(strat\_plan)}$ | $\Delta V$ | $t$ | $V_{(basic\ mode)}$ | $V_{(strat\_plan)}$ | $\Delta V$ |
|---|---|---|---|---|---|---|---|
| 1 | 87.3361 | 86.3361 | -1.0000 | 30 | 96.3218 | 95.3218 | -1.0000 |
| 2 | 70.9440 | 69.9439 | -1.0001 | 31 | 105.1011 | 104.1020 | -0.9991 |
| 3 | 51.4324 | 47.4322 | -4.0002 | 32 | 98.6620 | 97.6634 | -0.9986 |
| 4 | 56.3529 | 52.3533 | -3.9996 | 33 | 82.1931 | 81.1925 | -1.0006 |
| 5 | 59.2634 | 58.2642 | -0.9992 | 34 | 76.2280 | 75.2241 | -1.0038 |
| 6 | 73.3888 | 73.3893 | 0.0005 | 35 | 68.5178 | 68.5158 | -0.0020 |
| 7 | 95.2454 | 95.2459 | 0.0005 | 36 | 60.5149 | 60.5134 | -0.0015 |
| 8 | 92.6417 | 92.6422 | 0.0004 | 37 | 53.1263 | 53.1275 | 0.0013 |
| 9 | 95.5266 | 95.5275 | 0.0009 | 38 | 61.6513 | 61.6518 | 0.0006 |
| 10 | 70.1725 | 71.1725 | 1.0000 | 39 | 53.5085 | 53.5087 | 0.0002 |
| 11 | 58.4166 | 61.4166 | 3.0000 | 40 | 51.8387 | 51.8420 | 0.0033 |
| 12 | 56.4765 | 59.4765 | 3.0000 | 41 | 72.0274 | 72.0271 | -0.0002 |
| 13 | 61.8773 | 64.8773 | 3.0000 | 42 | 93.0773 | 93.0771 | -0.0002 |
| 14 | 71.8743 | 74.8685 | 2.9942 | 43 | 99.2277 | 99.2276 | -0.0001 |
| 15 | 52.4501 | 55.4537 | 3.0036 | 44 | 115.7921 | 115.7833 | -0.0089 |
| 16 | 53.9159 | 56.9331 | 3.0172 | 45 | 110.1212 | 110.1150 | -0.0062 |
| 17 | 84.0648 | 84.0583 | -0.0065 | 46 | 103.6429 | 103.6377 | -0.0052 |
| 18 | 114.4254 | 114.4221 | -0.0033 | 47 | 88.2175 | 87.2174 | -1.0002 |
| 19 | 132.1967 | 132.1964 | -0.0003 | 48 | 69.2655 | 68.2651 | -1.0004 |
| 20 | 153.9033 | 153.9141 | 0.0108 | 49 | 55.1417 | 54.1424 | -0.9993 |
| 21 | 164.5378 | 164.5295 | -0.0083 | 50 | 63.5110 | 62.5119 | -0.9991 |
| 22 | 150.0201 | 150.9962 | 0.9761 | 51 | 50.0201 | 49.0191 | -1.0009 |
| 23 | 140.7430 | 144.7379 | 3.9948 | 52 | 61.5773 | 60.5763 | -1.0010 |
| 24 | 115.0949 | 119.0896 | 3.9947 | 53 | 60.3298 | 59.3258 | -1.0040 |
| 25 | 87.0156 | 91.0163 | 4.0007 | 54 | 147.3276 | 147.3281 | 0.0004 |
| 26 | 100.5965 | 104.6030 | 4.0065 | 55 | 158.4131 | 158.4138 | 0.0008 |
| 27 | 87.5916 | 91.5951 | 4.0035 | 56 | 156.8698 | 156.8698 | 0.0000 |
| 28 | 76.0417 | 79.0439 | 3.0022 | 57 | 167.8964 | 167.8964 | 0.0000 |
| 29 | 76.2642 | 76.2632 | -0.0010 | | | | |

A figure 1 shows the experiment result of estimating the control mode $V_i(t)$ through the strategic planning.

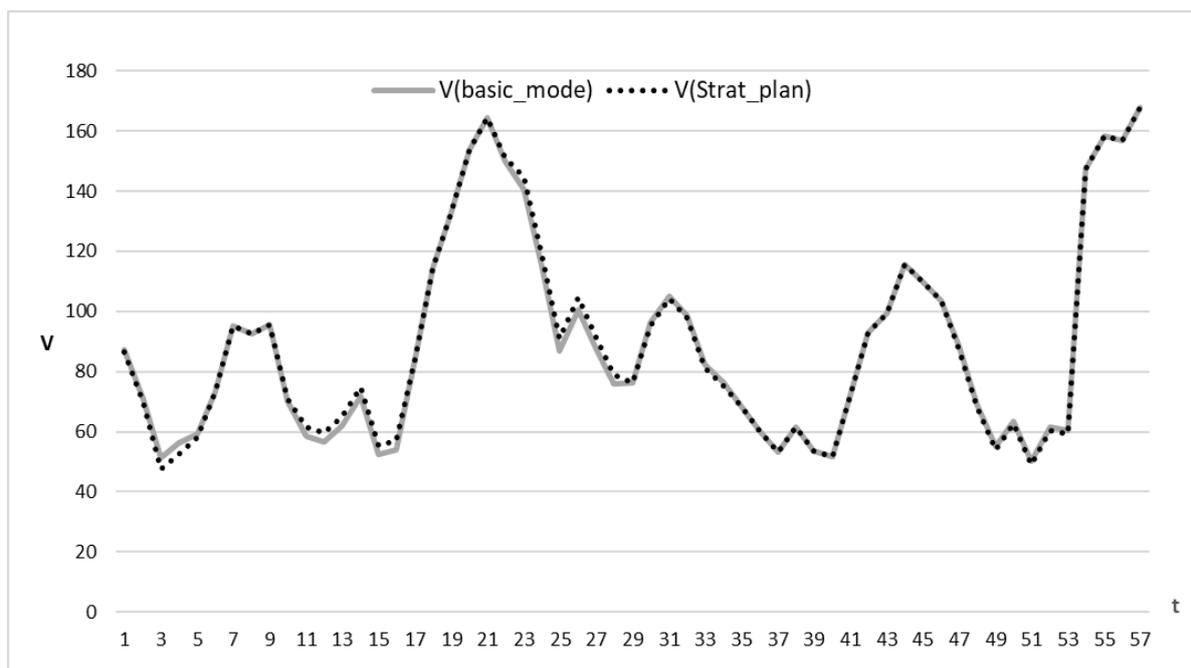

**Figure 1.** Indicator dynamics $V_i(t)$.

## 5. The discussion of the results

The enterprise had a simplified strategic planning system. Strategies have been identified and partially tied to the ongoing business processes of the enterprise. This simplified the work, since there was no need to create a separate center for strategic planning or additional training of company personnel. All managers and heads of departments were involved in the process of setting up the strategic planning. Their salaries were increased, additional resources were issued in the amount of 27,612 thousand rubles. Then the total costs of the enterprise for five years will amount to 5,669,054 thousand rubles.

The use of the integral indicator method allows taking into account the minimum changes occurring at the enterprise. For example, the difference between the implementation of control of the targets classifier and control of strategic planning.

## 6. Conclusion

The dynamics of the integral indicator reflects the implementation of management through strategic planning. In the standard mode of operation with the absence of linking strategies to the running business processes is $V_{(basic\_mode)}$ - 5,069.93 and with the linkage of strategies to business processes of the enterprise $V_{(strat\_plan)}$ - 5,089.90. Therefore, the assessment of the transition to object control through the strategic planning is estimated as $\Delta V = V^k_{strat\_plan} - V^k_{basic\_mode} = 19.96$. The purpose of the research has been achieved.